\documentstyle[12pt]{article}

  \newcommand{\const}{\rm const}
  \newcommand{\Var}{\rm Var}

  \newcommand{\Ent}{\rm  Ent}
   \newcommand{\rank}{\rm  rank}
   
\begin{document}

   \begin{center}

 {\bf  Exact exponential tail estimation for sums of independent } \\

 \vspace{3mm}

  {\bf centered random variables, under natural norming,} \par

 \vspace{3mm}

 {\sc with applications to the theory of  U-statistics }

\vspace{5mm}

 {\bf M.R.Formica,  E.Ostrovsky and L.Sirota.}

 \end{center}

\vspace{4mm}

 Universit\`{a} degli Studi di Napoli Parthenope, via Generale Parisi 13, Palazzo Pacanowsky, 80132,
Napoli, Italy. \\

e-mail: mara.formica@uniparthenope.it \\

\vspace{3mm}

Department of Mathematics and Statistics, Bar-Ilan University, \\
59200, Ramat Gan, Israel. \\

e-mail: eugostrovsky@list.ru\\

\vspace{3mm}

\ Department of Mathematics and Statistics, Bar-Ilan University, \\
59200, Ramat Gan, Israel. \\

e-mail: sirota3@bezeqint.net \\

\vspace{4mm}

\begin{center}

  {\bf Abstract} \par

\end{center}

\vspace{3mm}

 \ We derive in this short report the exact exponential decreasing tail of distribution for naturel normed
 sums of independent centered random variables (r.v.), applying the theory of Grand Lebesgue Spaces (GLS). \par
 \  We consider also some applications into the theory of $ \ U  \ $ statistics, where we deduce alike for the
 independent  variables   the refined
 exponential tail estimates for ones  under natural norming sequence.\par

\vspace{4mm}

\begin{center}

 { \sc Key words and phrases.}\par

\vspace{3mm}

 \end{center}

\vspace{3mm}

 \hspace{3mm} Probability, random variables (r.v.)  expectation, variance, distribution and tail of distribution,  $ \ U - \ $
statistic, independence, identical distributional, tail function, theorem of Fenchel - Moreau, convexity, generating function,
martingale representation, rang, exponential decreasing  tail estimations, examples,  sums of independent centered variables,
Lebesgue - Riesz and Grand Lebesgue spaces (GLS) and norms, ordinary and moment generating function, ordinary and exponential
estimation, Young - Fenchel transform, Young's inequality, confidence region, convergence, Laplace transform, kernel. \par

 \vspace{4mm}

 \section{Statement of problem. Previous results. Notations and conditions.}

 \vspace{4mm}

 \hspace{3mm} Let $ \ (\Omega = \{\omega\}, \ \cal{B}, {\bf P})   \ $  be probability space, with expectation \ {\bf E} \ and variance \Var.
Denote for arbitrary r.v. $ \ \eta \ $ its {\it tail function }

$$
T[\eta](t) \stackrel{def}{=} {\bf P}(|\eta| > t), \hspace{3mm} t \in (0,\infty).
$$

\vspace{3mm}

 \ Let also $ \   \{\xi_i\}, \ \xi = \xi_1, \ i = 1,2,\ldots,b \ $ be a sequence of {\it centered, identical distributed, with finite  non - zero variance
 and independent }random variables (r.v.)  Put as ordinary

$$
S_n = n^{-1/2} \sum_{i=1}^n \xi_i,
$$
and correspondingly

$$
Q[S_n](t) \stackrel{def}{=} T[S_n](t)  = {\bf P} (|S_n| > t), \  t \in (0,\infty).
$$

 \ It is known, see \cite{Ermakov etc. 1986}, \cite{MRFormica}  that if

\vspace{3mm}

\begin{equation} \label{m cond}
 \ \forall t > 1 \ \Rightarrow \ T[\xi](t) \le \exp \left(-t^m \right), \ \exists m = \const > 0,
\end{equation}

\vspace{3mm}

then

\vspace{3mm}

\begin{equation} \label{one dim estim}
T[S_n](t) \le \exp \left( \ - c_m \ t^{\min(m,2)} \ \right),  \ t \ge 1, \ c_m > 0.
\end{equation}

\vspace{3mm}

 \ Moreover, the last estimate is essentially non improvable. \par

 \vspace{3mm}

  \hspace{3mm} {\bf Our target in this report is to generalize the exponential estimate }  \ (\ref{one dim estim}) \   {\bf on the case of arbitrary
exponential  tail of distribution  behavior of source random variable as well as into the multidimensional case, more precisely, into the U - statistics. } \par

\vspace{3mm}

 \ We will rely essentially on the methods offered in the recent work \cite{Bakhshizadeh}; as well as on the theory of the so - called Grand Lebesgue Spaces
 (GLS), represented e.g. at the works \cite{anatriellofiojmaa2015}, \cite{anatrielloformicaricmat2016}, \cite{Buld Koz AMS}, \cite{Buldygin-Mushtary-Ostrovsky-Pushalsky},
\cite{Ermakov etc. 1986}, \cite{Fiorenza2000}, \cite{fiokarazanalanwen2004}, \cite{Fiorenza-Formica-Gogatishvili-DEA2018}, \cite{fioformicarakodie2017}, \cite{formicagiovamjom2015},
\cite{MRFormica},  \cite{Iwaniec-Sbordone 1992}, \cite{Koz Os},
\cite{liflyandostrovskysirotaturkish2010}, \cite{Ostrovsky1999}, \cite{Ostrovsky HIAT}, \cite{Ostrov Prokhorov}, \cite{OsSirUstat},
\cite{Samko-Umarkhadzhiev}, \cite{Samko-Umarkhadzhiev-addendum}.\par

 \vspace{4mm}

 \section{Refined tail estimations for ordinary sums.}

 \vspace{4mm}

 \ {\sc Given: } the (centered) random variable e.g. first in the our list $ \ \xi = \xi_1 \ $ is such that

 \vspace{3mm}

 \begin{equation} \label{g condition}
T[\xi](t) \le \exp(-g(t)), \ t \ge 0,
 \end{equation}

 \vspace{3mm}
where the function $ \ g = g(t)  \ $ is strictly convex positive continuous and such that $ \ \lim_{t \to \infty} g(t)/t = \infty  \  $ and

$$
g(0) = g'(0) = 0, \hspace{3mm} g''(0) \in (0,\infty).
$$

\vspace{3mm}

 \hspace{3mm} It is proved in particular in
 \cite{Koz Os},  \cite{liflyandostrovskysirotaturkish2010}, \cite{Ostrovsky1999}, that   $ \  \forall \lambda \in R \ \Rightarrow \ $

\vspace{3mm}

\begin{equation} \label{MGF}
{\bf E}\exp(\lambda \ \xi) \le \exp \left( \  g^*(C_1 \lambda) \ \right), \ C_1 = \const \in (0,\infty).
\end{equation}

\vspace{3mm}

 \ Here and henceforth $ \ g^*(\cdot) \ $ will be denote the famous Young - Fenchel, or Legendre transform for the  {\it arbitrary} function $ \ g(\cdot): \ $

 $$
 g^*(\lambda) \stackrel{def}{=} \sup_{t \ge 0} (|\lambda| t  - g(t)).
 $$

\vspace{3mm}

 \ We have for all the values $ \ \lambda \in R \ $

$$
{\bf E} \exp(\lambda \ S_n) = {\bf E} \exp \left( \ \lambda \ n^{-1/2} \ \sum_{i=1}^n \xi_i  \  \right) =
$$

$$
\prod_{i=1}^n {\bf E} \exp \left( \ \lambda \ n^{-1/2} \ \xi_i  \ \right) \le  \prod_{i=1}^n  \exp \left( \  g^*(C_1 \lambda/\sqrt{n}) \ \right) =  \exp \left[ \ \nu_n(\lambda)   \ \right],
$$
\
where
$$
\nu_n(\lambda) \stackrel{def}{=} \ n \ g^*(C_1 \ \lambda /\sqrt{n}), \ \lambda \in R.
$$

 \vspace{3mm}

 \ We will use the following modification of the famous Chernoff's  inequality, see \cite{Koz Os},
\cite{liflyandostrovskysirotaturkish2010}, \cite{Ostrovsky1999}, \cite{Ostrovsky HIAT}, \cite{Ostrov Prokhorov}:

\vspace{3mm}

\begin{equation} \label{Chernov}
{\bf P}(S_n > t) \le \exp \left(-\nu_n^*(t) \ \right), \ t > 0.
\end{equation}

\vspace{3mm}

 \hspace{3mm} Notice that

$$
\nu_n^*(t) = \sup_{z > 0} \left\{ \ \lambda t - n g^* \left( \ C_1 \frac{\lambda}{\sqrt{n}}   \ \right)   \ \right\} =
$$

$$
n \sup_{z > 0} \left\{ \ \frac{t z}{C_1 \ \sqrt{n}}  -  g^*(z) \ \right\} = n g^{**} \left( \ \frac{t}{C_1 \ \sqrt{n}}  \ \right).
$$

 \vspace{3mm}

 \   Theorem of Fenchel - Moreau says that (under our conditions) $ \ g^{**}(z) = g(z), \ $  therefore

\vspace{3mm}

\begin{equation} \label{Tail}
{\bf P}(S_n > t) \le \exp \left[ \ - n g \left( \ \frac{t}{C_1 \ \sqrt{n}}   \ \right)    \ \right].
\end{equation}

\vspace{3mm}

 \ To summarize.\par

 \ {\bf Theorem  1.1.} \ We conclude under formulated conditions

 \vspace{3mm}

\begin{equation} \label{One dim tail}
 \max \left\{ {\bf P}(S_n > t),   {\bf P}(S_n < - t) \ \right\} \le \exp \left[ \ - n g \left( \ \frac{t}{C_1 \ \sqrt{n}} \ \right) \ \right], \ t > 0,
\end{equation}

 \vspace{3mm}

and correspondingly

\vspace{3mm}

\begin{equation} \label{Bilateral}
Q[S_n](t) \le 2 \exp \left[ \ - n g \left( \ \frac{t}{C_1 \ \sqrt{n}} \ \right)  \ \right], \ t > 0.
\end{equation}

\vspace{3mm}

\hspace{3mm}

 \ {\bf Remark 1.1.} It is easily to verify that the mentioned above estimates \ (\ref{m cond}), (\ref{one dim estim})
follows immediately from   (\ref{Bilateral}). Thus, the  obtained now relation  (\ref{Bilateral}) is essentially non - improvable.\par

 \vspace{4mm}

 \section{Refined  exponential tail distribution estimations for U - statistics.}

 \vspace{4mm}

 \hspace{3mm} Let $ \  (\Omega,  \cal{B}, {\bf P} ) \ $ be again probabilistic space, which will be presumed sufficiently rich
when we construct examples (counterexamples). Let $ \ \{X(i)\}, i = 1, 2, . . . , n \ $
be independent identically distributed (i., i.d.) random variables (r.v.) with values in the
certain measurable space  $ \ (X, S), \ h = h(x(1), x(2), . . . , x(m)) \ $ be a {\it symmetric} measurable {\it centered} non - trivial
numerical function ({\it kernel} ) of $ \  m \ $  variables: $ \ h : X^m \to R, \ $

$$
{\bf E}h = {\bf E} h(X(1),X(2), \ldots,X(m)) = 0.
$$

\vspace{3mm}

 \ Introduce also as ordinary the variables

\vspace{3mm}

\begin{equation} \label{U stat def}
U_n = U(n, h, d) = U(n, h, d; \{X(i)\}) =
\end{equation}

\begin{equation} \label{U stat2}
\left( \stackrel{n}{m} \right)^{-1} \cdot \sum \ \sum \ldots \sum_{1 \le i(1) < i(2) \ldots i(m) \le n} h(X(i(1)), X(i(2)), \ldots, X(i(m))
\end{equation}

\vspace{3mm}
be a so-called $ \ U - \ $ statistic;

$$
 \deg h  = \deg U = m, \  \sigma^2(n) := \Var(U_n) \asymp n^{-r}, \ r = \rank(U),
$$

$$
\eta :=   h( X(1), X(2), \ldots, X(m)), \hspace{3mm} \beta^2 := \Var(\eta).
$$

$$
 k = k(m,n) := \Ent[n/m];
$$

where $ \ \Ent(Y) \ $ denotes the integer part of the variable $ \ Y; \ $

$$
T([U(n)],t)  \stackrel{def}{=} {\bf P} \left[(U_n - {\bf E} U_n)/\sigma(n) > t \ \right], \ t > 0.
$$

\vspace{3mm}

 \ We will use the following very important estimate which is grounded in the article  \cite{Bakhshizadeh}:

 \vspace{3mm}

 \begin{equation} \label{Bakhshizadeh}
 {\bf E} \exp (\lambda U_n) \le \left\{ \ {\bf E} \left[ \ \frac{\lambda \eta}{k} \right]  \ \right\}^k.
 \end{equation}

\vspace{3mm}

 \ Suppose as in the first section that

\vspace{3mm}

\begin{equation} \label{tail for the h}
T[\eta](t) \le \exp \left( \ - l(t)  \ \right), \ t \ge 0,
\end{equation}

\vspace{3mm}

where  the function $ \  l = l(t) \ $ obeys at the same properties as the introduced before the function $ \ g(\cdot). \ $
Therefore

\vspace{3mm}

\begin{equation} \label{MGF2}
{\bf E}\exp(\lambda \ \eta) \le \exp \left( \  l^*(C_2 \lambda) \ \right), \ C_2 = \const \in (0,\infty), \ \lambda \in R,
\end{equation}

\vspace{3mm}
and the inverse conclusion is true. \par
 \ Ww conclude substituting into   (\ref{Bakhshizadeh})

$$
{\bf E} \exp \left(\mu \ \sqrt{n} \ U_n \right) \le \exp \left\{ \ n \ l^* \left(\ \frac{\mu \ C_3(m)}{\sqrt{n}} \ \right) \ \right\}, \ \mu \in R.
$$

\vspace{3mm}

 \ We obtain finally quite analogously the proof of theorem  1.1 \par

\vspace{3mm}

 \ {\bf Theorem  2.1.} \ We conclude under formulated before notations and  conditions for all the positive values $ \ t \ $

 \vspace{3mm}

\begin{equation} \label{General case}
 \max \left\{ {\bf P}(\sqrt{n} \ U_n > t),   {\bf P}( \sqrt{n} \ U_n < - t) \ \right\} \le \exp \left[ \ - n \ l \ \left( \ \frac{t}{C_3(m) \ \sqrt{n}} \ \right) \ \right],
\end{equation}

 \vspace{3mm}

and correspondingly

\vspace{3mm}

\begin{equation} \label{Bilateral gen case}
Q[\sqrt{n} \ U_n](t) \le 2 \exp \left[ \ - n \  l \ \left( \ \frac{t}{C_3(m) \ \sqrt{n}} \ \right)  \ \right], \ t > 0.
\end{equation}

 \vspace{3mm}
 
  \ {\bf Remark 3.1.}  The case $ \ m = 1 \ $ correspondent to the considered before the one - dimensional case. 
So, theorem 2.1 is the direct generalization of the classical one - dimensional estimates. \par

\vspace{3mm}

\section{Concluding remarks.} 

\vspace{3mm}

 \hspace{3mm} It is interest in our opinion to generalize obtained results on the non - symmetrical kernels,
 as well as onto the Banach space random variables.\par 

 \vspace{6mm}

\vspace{0.5cm} \emph{Acknowledgement.} {\footnotesize The first
author has been partially supported by the Gruppo Nazionale per
l'Analisi Matematica, la Probabilit\`a e le loro Applicazioni
(GNAMPA) of the Istituto Nazionale di Alta Matematica (INdAM) and by
Universit\`a degli Studi di Napoli Parthenope through the project
\lq\lq sostegno alla Ricerca individuale\rq\rq .\par

\end{document}